\documentclass{article}
\usepackage[utf8]{inputenc}
\usepackage{amssymb}
\usepackage{mathptmx}
\usepackage{amsmath}
\usepackage{tikz}
\usepackage{natbib}
\usepackage{bussproofs}

\title{Counterfactuals in Branching Time: The Weakest Solution}
\author{Daniil Khaitovich \footnote{International Laboratory for Logic, Linguistics and Formal Philosophy, HSE, Moscow. Email:dkhaytovich@hse.ru}}
\date{ }
\begin{document}
\maketitle
\section*{Introduction}
Counterfactuals are subjunctive conditionals with false antecedent: "If I were 5 meters tall, I would be the tallest human being alive", "If there were no gravity, it would be harder to walk down the streets", "If I choosed to bet tail, I would win" etc. If we try to provide a formal semantics for these statements, we will find out that counterfactual conditionals require intensional logic: the validity of expressions of the form "if it were $\phi$, then it would be $\psi$" cannot be defined just in terms of truth-values of $\phi$ and $\psi$. Counterfactuals with true consequent may be both true (2) or false  (1), just as counterfactuals with false consequents:

\begin{equation}
    \mbox{If my heart stopped, I would still be able to run (0 $>$ 1; 0)}
\end{equation}

\begin{equation}
    \mbox{If I were bald, I would still be able to run (0 $>$ 1; 1)}
\end{equation}

Moreover, strict implication is not a good candidate too: if we try to use formulas of the form $\Box (\phi \rightarrow \psi)$ of any normal modal logic, where $\Box$ is interpreted as necessity operator, \textit{strengthening the antecedent} will be valid, i.e. 

\[
\Box (\phi \rightarrow \psi) \vdash \Box((\phi \land \chi)\rightarrow \psi)
\]

But it is strange to think that from

\begin{center}
    
\textit{If you had cut your finger, you would not need any medical help}
\end{center}

it follows that

\begin{center}

\textit{If you had cut your finger and your throat, you would not need any medical help}
\end{center}

The most common semantic framework of logic for counterfactual conditionals, which escapes these pitfalls, is so called comparative similarity analysis, proposed by Robert Stalnaker \cite{stalnaker1968} and David Lewis \cite{lewis1973}. The basic idea may be expressed as follows:
\begin{center}

\textit{"If it were $\phi$, then it would be $\psi$"} is true at some possible world $w$ iff in all $\phi$-worlds accessible from $w$, which differ minimally from $w$ itself, $\psi$ is true. 
\end{center}

Since we evaluate the expression in regards of some subset of the accessible worlds, which depends on the antecedent, \textit{strengthening the antecedent} is easy to falsify. 

The proper formal definition of the comparative similarity semantics and  the corresponding logic is well-known established results.
But new questions arise when we think of adding temporality in the framework. Since the very idea of counterfactual statements might imply alternative possibilities (we will discuss it in details in Section 2), the most natural way to represent the notion is to define the logic of counterfactuals with respect to branching time structures. But then we face both formal and conceptual problems. In this paper we will give an overview of these concerns and propose the minimal solution. In Section 1, the simplest logic of counterfactuals \textbf{P} is presented. In Section 2, we discuss main metaphysical problems connected with time and counterfactual conditionals: determinism versus indeterminism, the difference between historical and plain counterfactuals. Section 3 presents Ockhamist branching time temporal logic  OBTL, which has several nice properties, suitable for counterfac-\-tuals. At the last Section we combine both logics presented in Section 1 and Section 3 in regards of concerns, observed in Section 2.  

\section{Minimal logic P}

We will precisely present the logic of counterfactuals and its comparative similarity semantics. The logic is the simplest system \textbf{P}. Other popular systems are usually \textbf{P}'s extensions: we can add more axioms depending on our metaphysical views on possible worlds and similarity relations on them. We will look at \textbf{P} for the sake of simplicity. 

The syntax of our language defined in BNF form as follows. For a countable set of propositions $Var$ and $p \in Var$:

\[\phi, \psi : = p \: | \: \neg \phi \: | \: \phi \lor \psi \: | \: \phi \rightsquigarrow \psi \]

$\phi \rightsquigarrow \psi$ means "if it were $\phi$, then it would be $\psi$". Other Boolean connectives are interpreted in a standard way. Semantically our logic is defined on the frames of the form:

\[
\mathcal{F}_1 = \langle W, \langle \lessdot_w \rangle_{w \in W}, f \rangle
\]

\begin{itemize}
    \item $W \neq \emptyset$ is a non-empty set of possible worlds.
    \item $\langle \lessdot_w \rangle_{w \in W}$ is a tuple of strict partial ordering on $W$, defined for every world $w \in W$. The orderings represent comparative similarity relations. For example, for some $x, y, w \in W$, $x \lessdot_w y$  means that $x$ is more similar to $w$ than $y$. For every $w \in W$, its corresponding comparative similarity relation is defined on a subset of worlds $W_w \subseteq W$, which are accessible from $w$. 
    \item $f: W \mapsto \{ \lessdot_w \}_{w \in W}$ is a function, which maps a corresponding strict partial ordering to every world.
\end{itemize}

As usual, a frame is extended to a model $\mathcal{M} = \langle\mathcal{F}_1 ,\nu \rangle$, where $\nu: Var \mapsto 2^W$ is an evaluation function, mapping each proposition to a subset of possible worlds, in which that proposition is true.  Formal semantics is natural here:
$\\$

\noindent$\mathcal{M}, w \models p  \mbox{ iff } w \in \nu(p)$\\
        $\mathcal{M}, w \models \neg \phi  \mbox{ iff } \mathcal{M}, w \not\models \phi $\\
        $\mathcal{M}, w \models \phi \lor \psi \mbox{ iff } \mathcal{M}, w \models \phi \mbox{ or }  \mathcal{M},w \models \psi$ \\
		$\mathcal{M}, w \models \phi \rightsquigarrow \psi \mbox{ iff } \mathcal{M}, v \models \psi \mbox{ for every } \lessdot_w\mbox{-maximal world } v\in W_w \mbox{, such that } \mathcal{M}, v \models \phi$
$\\$

The adequate logic of the frames defined above is namely \textbf{P}. It is sound and complete on the class of frames defined earlier \cite{burg1981}:
$ \\$

 \noindent(PL)  All tautologies of classic propositional logic  \\
                (CI)  $\phi \rightsquigarrow \phi$\\
                (CC) $((\phi \rightsquigarrow \psi) \land (\phi \rightsquigarrow \chi)) \rightarrow (\phi \rightsquigarrow (\psi \land \chi))$\\
                (CW)  $(\phi \rightsquigarrow \psi) \rightarrow (\phi \rightsquigarrow (\psi \lor \chi))$\\
                (SA) $((\phi \rightsquigarrow \psi) \land (\psi \rightsquigarrow \chi)) \rightarrow ((\phi \land \psi)\rightsquigarrow \chi) $\\
                (AD) $((\phi \rightsquigarrow \chi) \land (\psi \rightsquigarrow \chi)) \rightarrow ((\phi \lor \psi)\rightsquigarrow \chi)$ \\
                (REA) If $\vdash \phi \leftrightarrow \psi$ then $(\phi \rightsquigarrow \chi) \leftrightarrow (\psi \rightsquigarrow \chi)$\\
                (REC) If $\vdash \phi \leftrightarrow \psi$ then $(\chi \rightsquigarrow \phi) \leftrightarrow (\chi \rightsquigarrow \psi)$\\
\noindent (MP) \AxiomC{$\phi, \phi \rightarrow \psi$} \UnaryInfC{$\psi$} \DisplayProof
$\\$
                
As we have noted before, other axioms can be added, if we impose more metaphysically motivated restrictions on comparative similarity relation. But we will omit this discussion and work  with that simplest system. 
\section{Metaphysical concerns}

Now we need to apply the logic we have defined earlier to temporally sensitive expressions. But before going into details of temporal logic framework we want to use, there is a need to figure out several metaphysical concerns about time and alternative courses of events. 

\subsection{Plain versus Historical Counterfactuals}

Some subjunctive conditionals may be interpreted in two different fashions. For example, 

\begin{center}
    \textit{If I were 3 meters tall, I would be the tallest person in the world}
\end{center}

We can understand it as just the idea about the notion of being me: if we take that notion and change one of the properties it implies -- namely, my height -- that new notion of me necessarily will imply the property of being the tallest man in the world. The other way is to think of it as a thought about actual me and alternative course of events. If we imagine some alternative history, which is as similar to the real history as possible, where I am in fact 3 meters tall, in that alternative history I will be the tallest. But in that interpretation, it may be the case that the possible course of events resulting in me being that tall presuppose me getting other genetic heritage, consequently, it may presuppose other course of evolution of human beings, where 3 meters is an average height of a mature man, so the original statement might be false. 

Or another example:

\begin{center}
    \textit{If John were not late, he wouldn't drive when it's dark}
\end{center}

Again, we can understand it as  the very notion of John not being late presupposes that property of not driving in the darkness, or we can think of just the alternative course of events in the real world when John was not late, while every other aspect of the history stays the same as soon as it is compatible with John not being late. It seems obvious that the second interpretation is more natural in that context. In some other cases, the first way of interpreting coutnerfactuals will be more natural: 

\begin{center}
    \textit{If the square root of 2 were a rational number, it would be possible to express it as a ratio of 2 integers}
\end{center}

It is even hard to imagine that alternative course of real world evolution where the square root of 2 would be a rational number, so the expression is more likely to be a purely conceptual one. Wawer and Wroński proposed the difference between \textit{plain} and \textit{historical} counterfactuals to grasp that two interpretations \cite{waw2015}. If the pure comparative similarity analysis is enough for plain counterfactual, we need to compose it with a certain temporal logic somehow to grasp the historical one. 
\subsection{Local Miracles and Indeterminism}

The next problem is indeterminism associated with historical counterfactual conditionals: we are considering the possibility of alternative course of events in the actual world, i.e. at some point, given the  same past and settled laws of nature, the future might be different. 

Several authors argue that we should stay either neutral or taking the side of determinism when we are dealing with counterfactuals \cite{mul2007}. The most remarkable case is Lewis himself: in his article "Are We Free to Break the Laws?", Lewis argues for compatibilism (an idea that causal determinism is compatible with free will) with the help of counterfactual definition of ability to act otherwise \cite{lewis1981}. 

We will consider the next premises:

\begin{enumerate}
    \item The events of the past ($H$) and the laws of nature ($L$) causally determine future course of events
    \item Every human action, for example, raising one's hand ($R$), as an event, is predetermined by $H$ and $L$: $\Box((H \land L) \rightarrow R)$.
    \item Therefore, if one's acting against their predetermination, for example, not raising one's hand ($\neg R$), then either our past were different ($\neg H$), or the laws of nature were broken ($\neg L$).
\end{enumerate}

Lewis claiming that we do have a possibility to act otherwise, i.e. to do something, that if we did it, the laws of nature would be broken ($\neg R \rightsquigarrow \neg L$). But it does not mean that we are able to break the laws of nature: that hypothetical law-breaking event might occur before our action and not be caused by us.

As Placek and Müller noticed\cite{mul2007}, the position seems a bit ambivalent: while stating determinism, Lewis speaks of \textit{local miracles}, i.e. the events, after which the course of events starts to deviate from the path, predetermined by the past and laws of nature. This picture, where we have different histories and forks, is almost isomorphic to indeterministic theory we need for historical counterfactuals. We will see the similarities in the next section.

\section{Branching Time Temporal Logic}
                
As we have noticed before, the very notion of historical counterfactual statements implies certain kind of indeterminism, especially when it comes to time and so called future contingents. Luckily we have a well-developed formal theory of indeterminism without thin red line, i.e. the theory which states that alternative courses of events are real and there is no unique actual scenario that happened, happens and will continue happening. It is the theory of Ockhamist branching time, developed by Arthur Prior, Richard Thomanson, Nuel Belnap and other contributors \cite{belnap2001}. If we don't want to stick to indeterminism and save Lewisean view, we can regard this alternative courses of events not as real, but rather as mere possibilities of local miracles.

The main advantage of Ockhamist branching time theory in the context of counterfactuals is that it allows both expressions about time and historical possibility/necessity. "It will be hot tomorrow, but it is possible that it will be cold" is a consistent, perfectly well-formed statement, so that "if it were the case that it will be cold, we would be in need to find a warm coat" will be. 

The basic idea is to represent time  with the help of generalized McTaggart's B-series: we postulate a non-empty set of moments $M$, which is ordered by a precedence relation: $m_1 \succ m_2$ means that some moment $m_1$ comes earlier than another moment $m_2$. This relation should be partially ordered, i.e. it should allow some moments to be incompatible in regards to precedence: we will think of this pairs of moments in the next manner: if $m_1 \not\succ m_2$ and $m_2 \not\succ m_1$, then $m_1$ and $m_2$ belong to different \textit{histories}, i.e. full courses of events, how things might be. Then we may define some linearly ordered subsets of $M$, which will represent histories. Obviously, the subsets should be \textit{maximal} linearly ordered subsets of $M$ to illustrate the full histories. So that, the flow of time may be represented by a tree-like structure, where every node designates moments and edges -- precedence relations, every longest path -- a history (see Figure 1 for an illustration).

Another important feature is \textit{backward linearity}: for any $m_1, m_2, m_3 \in M$, if $m_1 \succ m_3$ and $m_2 \succ m_3$ then $m_2 \succ m_1$ or $m_1 \succ m_2$ or $m_1 = m_2$. It means that past is settled: alternative possibilities lies only in future.

Now we can formally define our  Ockhamist branching time temporal logic OBTL. The syntax can be presented in BNF form as follows. Given a countable set $Var$, for any propositional variable from that set $p \in Var$:

\[
\phi, \psi : = p \: | \: \neg \phi \: | \: \phi \lor \psi \: | \: G \phi \: | \: H \phi \: | \: \Box \phi  
\]

Atomic propositions and Boolean connectives have standard meaning. $G \phi$ reads as "at every moment in the future, $\phi$", $H\phi$ -- "at every moment in the past, $\phi$", $\Box \phi$ -- "it is historically necessary that $\phi$", which means that in all possible alternative histories, it is $\phi$ at the moment. As usual, all modal operators have duals:

\begin{center}
    $G \phi \equiv \neg F \neg \phi $\\
    $H \phi \equiv \neg P \neg \phi $\\
    $\Box \phi \equiv \neg \Diamond \neg \phi$
\end{center}

$F\phi$ stands for the truth of $\phi$ at \textit{some} moment in the future, $P \phi$ -- $\phi$ at some moment in the past and $\Diamond \phi$ -- for $\phi$ is true at the moment in at least some possible history.

So that, the Ockhamist branching time temporal logic is defined on a frame \[\mathcal{F}_2 = \langle M, H, \succ \rangle\] where $M = \{ m_1, m_2,... \}$ is a non-empty countable set of moments. $\succ \subseteq M^2$ is a partial ordering defined on $M$. A set of histories is $H = \{ h_1, h_2,...\}$, where each history is a maximal linearly $\succ$-ordered subset of $M$. A model is  $\mathcal{M} = \langle \mathcal{F}_2, \nu \rangle$, where $\mathcal{F}$ is a frame and $\nu: Var \mapsto 2^{M \times H}$ is a standard evaluation function, mapping a set of moment/history pairs (we will call them \textit{points}) to every atomic proposition. The key difference with linear temporal logics is that evaluation depends not only on moments, but on histories as well: the same moment may satisfy different formulas depending on what possible history is taking place. 
\begin{figure}[h]
    \centering
    \includegraphics{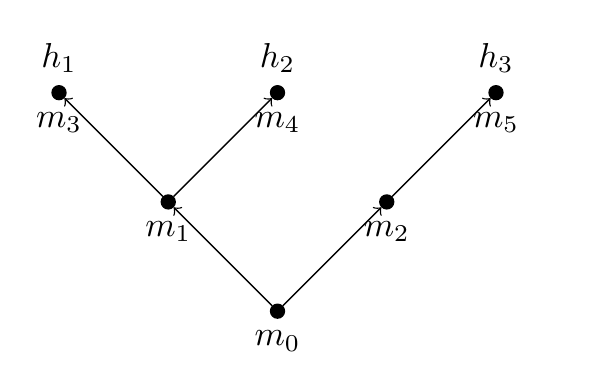}
    \caption{An example of BT frame}
    \label{fig:my_label}
\end{figure}

By $H_m$ we will denote the set of histories, running through the moment $m$: $H_m = \{ h \: | \: m \in h \}$. Now we can provide semantics for OBTL language.
$\\$

 \begin{align*}
    \noindent \mathcal{M}, m/h &\models p \mbox{ iff } m/h \in \nu(p) \\
     \mathcal{M}, m/h &\models \phi \lor \psi \mbox{ iff } \mathcal{M}, m/h  \models \phi \mbox{ or } \mathcal{M}, m/h  \models \psi\\
     \mathcal{M}, m/h &\models \neg \phi \mbox{ iff } \mathcal{M}, m/h  \not\models \phi\\
     \mathcal{M}, m/h &\models P\phi \mbox{ iff } \exists m'(m'\succ m \land \mathcal{M}, m'/h \models \phi) \\
      \mathcal{M}, m/h &\models F\phi \mbox{ iff } \exists m' (m \succ m' \land \mathcal{M}, m'/h \models \phi)\\
      \mathcal{M}, m/h &\models \Box \phi \mbox{ iff } \forall h' \in H_m (\mathcal{M}, m/h' \models p)\\
 \end{align*}

 When it comes to complete axiomatization for OBTL, it is still an open issue. We will use the logic proposed by Mark Reynolds in \cite{rey2002}. Moreover, we are going to treat it as a black box: no interference in OBTL axioms will be made, we will just extend the axiomatics with other systems and axioms, despite what the OBTL axiomatics is actually like. We will stick to the next variant proposed by Reynolds, which is proved to be sound:
$\\$

\noindent  All tautologies of classic propositional logic  \\
K4 for both $G$ and $H$ operators\\
S5 for $\Box$ operator\\
(L1) $\phi \rightarrow GP\phi$\\
(L1') $\phi \rightarrow HF \phi$\\
(L2) $F \phi \rightarrow G(F\phi \lor \phi \lor P \phi)$\\
(L2') $P \phi \rightarrow H(P \phi \lor \phi \lor F \phi)$\\
(L3) $\Box H \phi \equiv H \Box \phi$\\
(L4) $P \Box \phi \rightarrow \Box P \phi$\\
(L5) $\Box G \phi \rightarrow G \Box \phi$\\\
(L6) $G \bot \rightarrow \Box G \bot$
$\\$
The logic is closed under the next inference rules:
$\\$

\noindent(MP)
\AxiomC{$\phi, \phi \rightarrow \psi$}
\UnaryInfC{$\psi$}
\DisplayProof
\\
(RN)
\AxiomC{$\phi$}
\UnaryInfC{$\Box \phi$}
\DisplayProof
\AxiomC{$\phi$}
\UnaryInfC{$G \phi$}
\DisplayProof
\AxiomC{$\phi$}
\UnaryInfC{$H \phi$}
\DisplayProof

\vspace{3mm}
\noindent(IRR)
\AxiomC{$(p \land H \neg p)\rightarrow \phi$}
\UnaryInfC{$\phi$}
\DisplayProof
 where $\phi$ does not contain $p$.
 
 \section{Combining logics}
 \subsection{Degrees of strictness}
 Now we face our main task: how to combine \textbf{P} with OBTL in the most meaningful and formally correct way? The biggest conceptual concern is the scope of comparative similarity relations: which moment/history pairs should be comparable? In the original theory proposed by Lewis, we can compare logically/nomically accessible worlds (depending on how strict we want our counterfactual to be), i.e. we can compare how arbitrary two worlds $w_1, w_2$ are similar to $w$ if $w_1$ and $w_2$ are logically consistent/share the same laws of nature with $w$ \cite{lewis1973}. 
 
 When it comes to historical counterfactuals, there are plenty of candidates to be the adequate  restriction. The first one is \textit{co-presence}: we can compare only those moment/history pairs, which occur at the same time in different histories. For example, if we want to evaluate the sentence \textit{"If I were polite today, you would not complain about how annoying I am with being rude all the time"}, we need to bear in mind exactly that points of alternative histories, where 1) I was polite at the given day and 2) which are co-present with the moment of utterance of the original sentence: if we consider the history, where I am polite today, it is still hard to tell if I would face the complaint or not in the nearest future or 3 years ago in that alternative history. We can even strengthen the co-presence restriction to historical accessibility, stating than only historically accessible moment/history pairs are comparable.

\begin{figure}[h]
    \centering
    \includegraphics[scale=0.7]{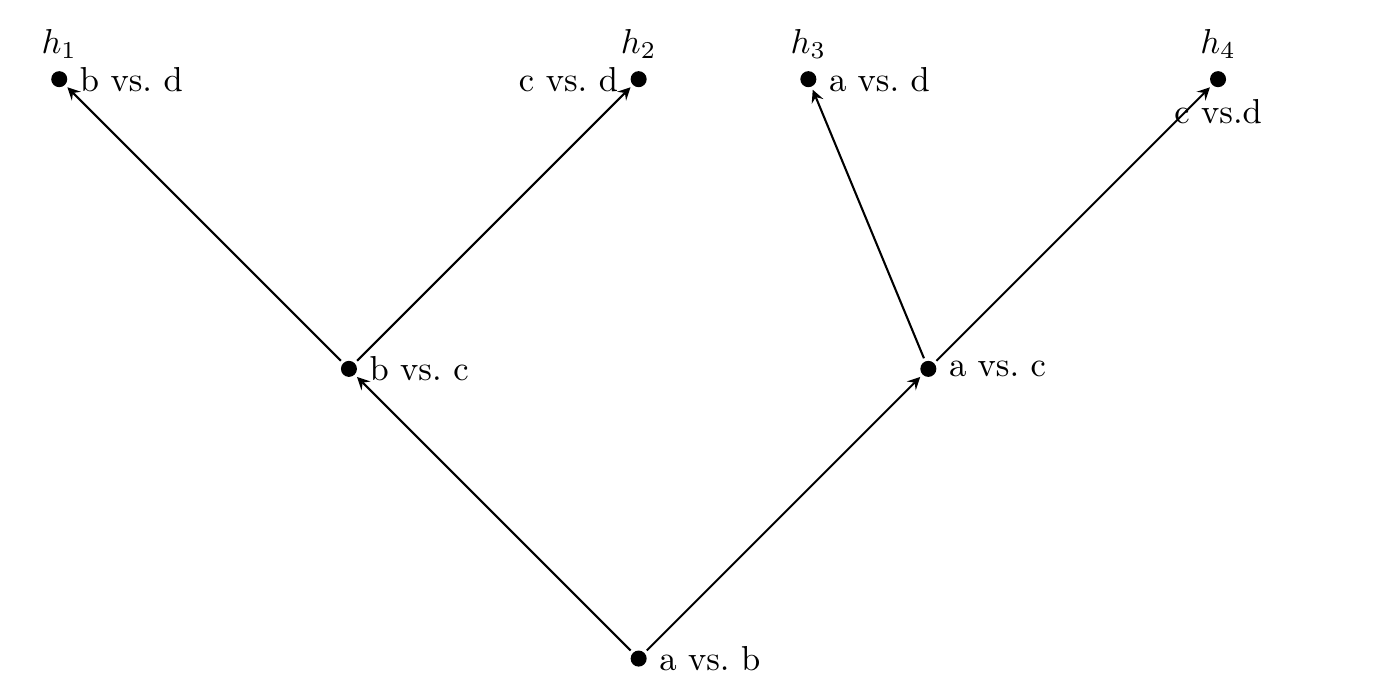}
    \caption{Tournament tree}
    \label{fig:my_label}
\end{figure}

 To illustrate the distinction, let us consider the next example. We have a tournament of philosophers: Aristotle (a), Berkeley (b), Chrysippus (c) and Descartes (d) are participating. On the Figure 2, we can see an Ockhamist branching time model of the tournament. It starts with the battle of Aristotle and Berkeley, the winner is going to compete with Chrysippus in the next round and the defeater is going to challenge Descartes at the final stage (for the matter of simplicity, we haven't graphically  illustrated the evaluation function of that model, but the truth/falsity of the propositions we are going to mention should be obvious: for example, at the point $(b \: vs. \: c/h_1)$ it is true that Berkeley has won the previous round and will win the current one). For example, Aristotle wins the first match and faces Chrysippus. The fan of stoicism states:
 \begin{center}

 \textit{"If Chrysippus were debating Berkeley, he would easily refute his claims and win"}. 
 \end{center}
 
 If we take co-presence as the only restriction for comparative similarity relation, then we need to take that moment/history pairs, which are the most similar to $(a \: vs.\:  c/h_3)$ and in which Chrysippus confronts Berkeley, from the next set:
 
 \[\{(a \: vs.\:  c/h_3), (a \: vs.\:  c/h_4), (b \: vs.\:  c/h_1), (b \: vs.\:  c/h_2) \}\]
 
 Basically, Chrysippus supporter is right with his counterfactual statement, if $(b \: vs.\:  c/h_2)$ is more similar to $(a \: vs.\:  c/h_3)$ than $(b \: vs.\:  c/h_1)$. 
 
 But if we want to compare only historically accessible points, then we need to consider only $\{(a \: vs.\:  c/h_3), (a \: vs.\:  c/h_4) \}$. None of these points models the situation of Chrysippus debating Berkeley, so the given counterfactual is vacuously true: it is historically impossible for the battle mentioned in antecedent to occur. In that case, the Chrysippus supporter should be more accurate and reformulate his claim:
 
 \begin{center}

 \textit{"In the past it was true that if Berkeley had defeated Aristotle and were debating Chrysippus, the latter one  would easily refute Berkeley's claims and win"}. 
 \end{center}
 
 Intuitively, the reformulation is more precise: by saying \textit{"If Chrysippus were debating Berkeley..."}, it was meant the case in which the previous round, Aristotle versus Berkeley, ended up differently.
 
 Obviously, historical accessibility is a stronger notion than co-presence: all historically accessible points are co-present. But even co-presence is not so simple. Let us consider the next counterfactual:
 
 \begin{center}
     \textit{"If I were 17 years old, I would not be able to buy a bottle of wine"}.
 \end{center}

 In terms of comparative similarity analysis, if we take co-presence as a valid restriction, we shall interpret that sentence in a very unnatural way: in some history, as similar to the actual one as possible, where time goes in a deviant way and I am younger than actual me \textit{at the same moment}, I am not able to purchase a bottle of wine. While the most intuitively appealing interpretation will be that simple: when I was 17, I couldn't buy a bottle of wine. 
 
 \begin{figure}[h]
     \centering
     \includegraphics{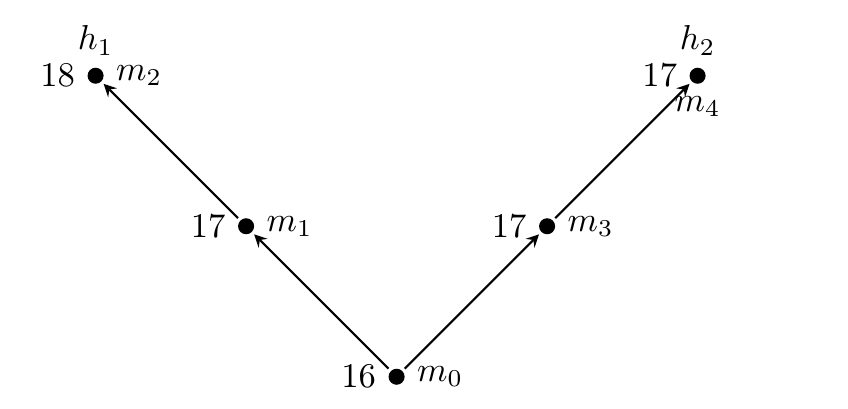}
     \caption{The bottle of wine example}
     \label{fig:my_label}
 \end{figure}
 
 Figure 3 illustrates the corresponding model. If we evaluate our counterfactual at $(m_2/h_1)$, where I am 18 years old, than  we should look either at co-present $(m_4/h_2)$, where I am 17 years old somehow (maybe I have spent some time on the other planet of Solar system after 17), or at the $(m_1/h_1)$, the actual past when I was 17. The example is not so artificial: it is common to express statements about actual past in the counterfactual manner. \textit{"If I were younger, I would be able to...", "If today were yesterday, he might not be late"}, etc. But not all "if today were yesterday" cases are inconsistent with co-presence: sometimes it is important to save some essential features of the current moment. \textit{"If he were younger, he would marry her"} may be an example: obviously, it is not synonymous to \textit{"He would marry her in the past, when he was young"}. 
 
 To conclude, historical counterfactuals are extremely vague, so we have to be accurate not to impose way too strict definitions. That is our motivation to favour for as weak definition as possible, allowing different context-dependent evaluations. We will regard  neither  co-presence nor historical accessibility as a necessary restriction to comparative similarity relations. Stronger system with co-presence was proposed by Canavotto: \cite{canavotto2020}.
 
 Another motivation for us to abandon the restrictions is a metaphysics behind Ockhamist branching time theory: all moment/history pairs occur in the only one world, \textit{our world}. Regarding Lewis theory, all points should be both logically and nomically accessible for each other. 
 \subsection{\textbf{P}+OBTL}
 
 Now we will formally describe the logic $\mathcal{L}_{PBT}$, corresponding to our understanding of historical counterfactuals. The language of the logic may be defined as follows. For a countable set of propositions $Var$ and $p \in Var$:
 
 \[
 \phi, \psi := p \: | \: \neg \phi \: |\: \phi \lor \psi \: | \: \phi \rightsquigarrow \psi \: | \: G \phi \: | \: H \phi \: | \: \Box \phi  
 \]
 
  As we can see, it is just a combination of \textbf{P} and OBTL languages. The $\mathcal{L}_{PBT}$ frame is a tuple:
  
  \[
  \mathcal{F}_3 = \langle M, H, \succ, \langle \lessdot_{m/h}\rangle_{m/h \in I}, f  \rangle
  \]

\begin{itemize}
    \item $M \neq \emptyset$ is a non-empty countable set of moments, $\succ$ is a partial ordering on $M$, $H\subseteq 2^M$ is a set of histories, i.e. maximal linearly $\succ$-ordered subsets of $M$, just like in OBTL frame. We will denote a set of all points in the frame as $I \subseteq (M \times H)$. 
    \item $\langle \lessdot_{m/h}\rangle_{m/h \in I}$ is a tuple of strict partial orderings on $I$, defined for every point $m/h \in I$. The orderings represent comparative similarity relations.
\end{itemize}

As usual, a frame $\mathcal{F}_3$ is extended to a model $\mathcal{M}= \langle \mathcal{F}_3, \nu \rangle$ by adding a valuation function $\nu: Var \mapsto 2^{M \times H}$, which maps every proposition to a set of points, in which the proposition is true.

Semantics stays the same for propositions, Boolean connectives and modal operators, as it was in OBTL.  Naturally, semantics for $\rightsquigarrow$ is defined as follows:
$\\$

\noindent $\mathcal{M}, m/h \models \phi \rightsquigarrow \psi \mbox{ iff } \mathcal{M}, m'/h' \models \psi$ for every $\lessdot_{m/h}$-maximal point $m'/h'$, such that \\ $\mathcal{M}, m'/h' \models \phi$
$\\$

We can see that semantically $\mathcal{L}_{PBT}$  is nothing more than a fusion of minimal counterfactual logic \textbf{P} and Ockhamist branching time temporal logic OBTL. Consequently, if we combine \textbf{P} axioms with OBTL axioms (i.e. join two sets of axioms and define it as closed under the rules of inference of both systems), we will preserve completeness \cite{fine1996}.

In order to show it explicitly, we need to redefine our frame in the next manner:

\[
\mathcal{F}_{3'} = \langle W, \succ, R_{\Box}, \langle \lessdot_w \rangle_{w \in W}, f \rangle
\]

\begin{itemize}
    \item $W \neq \emptyset$ is non-empty countable set of possible worlds or states
    \item $\succ \subseteq W \times W$ is a partial ordering on $W$
    \item $R_{\Box} \subseteq W \times W$ is an equivalence (i.e. transitive, reflexive and symmetric) relation, defined on $W$: for arbitrary $w_1, w_2 \in W$, if $w_1 R_{\Box}w_2$, then $w_2$ is historically accessible from $w_1$, i.e. $w_1$ and $w_2$ are points sharing the same moment, but contained in different histories
    \item $ \langle \lessdot_w \rangle_{w \in W}$ is a tuple of strict partial orderings defined on $W$ for every possible state. Their definition and meaning stays the same as in \textbf{P}.
    
\end{itemize}

Now we can clearly see that the frame $\mathcal{F}_{3'}$ shares the same domain $W$ and the same tuple of relations $\langle \lessdot_w \rangle_{w \in W}$, as it was in $\mathcal{F}_1$ case. We have added $R_{\Box}$ and $\succ$ relations, corresponding to OBTL frame $\mathcal{F}_2$: the only difference is that instead of defining histories as $\succ$-maximal subsets of $W$ and evaluating expressions on points as moment/history pairs, we have defined an equivalence relation on points, sharing the same moment. It won't affect our logic, since we can present a family of sets of moment/history pairs sharing the same moment $ \{\{ m/h \: | \: h \in H_m \} \} _{m \in M}$ as a partition of the set of all points $W$, and it is trivial to show that we can define an equivalence relation on $W$, corresponding to that partition, and it is exactly $R_{\Box}$. 

After extending the frame $\mathcal{F}_{3'}$ to a model in a usual manner, the only semantically different definition we will see is that for $\Box \phi$ formulas. Other definitions stay the same with the only change is that we use states $w \in W$ instead of moment/history pairs $ m/h \in (M \times H )$:

\begin{equation*}
 \mathcal{M}, w \models \Box \phi \mbox{ iff } \forall w' \in W (w R_{\Box} w' \rightarrow \mathcal{M}, w' \models \phi)  
\end{equation*}

Now we can observe that:
\begin{enumerate}
    \item $\mathcal{L}_{PBT}$ language is freely generated by the union of signatures $\textbf{P}$ and $OBTL$
    \item $\mathcal{L}_{PBT}$ model is obtained by combining relations defined on the same domain $W$, saving the same evaluation function $\nu$
\end{enumerate}

So that, by a results proved in \cite{fine1996}, the logic $\mathcal{L}_{PBT} = \textbf{P} \bigoplus OBTL$, which is generated by the simple union of axiom schemas and rules of inference of both $\textbf{P}$ and $OBTL$, is sound and complete with regards to $\mathcal{F}_{3'}$ frames (as soon as $OBTL$ is complete with regards to $\mathcal{F}_2$ frames).

No further constraints on relations interplay were imposed, hence, no multimodal axioms is needed. $\mathcal{L}_{PBT}$ is a conservative extension of both $\textbf{P}$ and $OBTL$.

\section{Conclusion} 

We have presented a formal analysis of temporally sensitive counterfactual conditionals, using the fusion of Ockhamist branching time temporal logic and minimal counterfactual logic \textbf{P}. We have presented an overview of the problems, both formal and metaphysical, which occur in combination of counterfactuals with temporal logics, and showed our motivation to use exactly that formal apparatus. 

Nevertheless, we are aware of the variety of problems we have refrained from touching. We haven't elaborated on formal constraints on comparative similarity relations and the corresponding axioms, as well as philosophical motivations to impose them. The enormously big discussion on determinism versus indeterminism was slightly mentioned as well.  Moreover, we haven't paid any attention to spatial aspect of the problem, observing only temporal one. We have knowingly left this white spaces in our overview to concentrate on a single aspect. This obstacle leaves the possibility to continue our work.

%clearpage
 %\textbf{Any opinions or claims contained in this Working Paper do not necessarily reflect the views of National Research University Higher School of Economics.}

 %textbf{© Khaitovich, 2021}

\end{document}